\documentclass[11pt,a4paper]{article}

\usepackage[utf8]{inputenc}
\usepackage[T1]{fontenc}
\usepackage[english]{babel}

\usepackage{amsmath,amssymb,amsthm,booktabs}
\usepackage[margin=2.5cm]{geometry}

\usepackage{tikz}
\usetikzlibrary{arrows.meta}

\usepackage{todonotes}
\usepackage[colorlinks=true,
linkcolor=blue,
citecolor=blue,
urlcolor=blue,
psdextra]{hyperref}

\newtheorem{theoreme}{Theorem}
\newtheorem{proposition}[theoreme]{Proposition}
\newtheorem{lemme}[theoreme]{Lemma}
\newtheorem{corollaire}[theoreme]{Corollary}
\newtheorem{definition}[theoreme]{Definition}
\newtheorem{exemple}[theoreme]{Example}
\theoremstyle{remark}
\newtheorem{remarque}[theoreme]{Remark}
\newtheorem*{notation_env}{Notation}

\newcommand{\transpose}[1]{{}^t #1}
\newcommand{\N}{\mathbb{N}}
\newcommand{\Z}{\mathbb{Z}}
\newcommand{\R}{\mathbb{R}}
\newcommand{\C}{\mathbb{C}}
\newcommand{\som}{\operatorname{sum}}
\newcommand{\parent}{\operatorname{par}}
\newcommand{\bl}{\boldsymbol{\ell}}
\newcommand{\bs}{\mathbf{s}}
\newcommand{\ee}{\mathbf{e}}
\newcommand{\w}{\mathbf{w}}
\newcommand{\bb}{\mathbf{b}}
\newcommand{\rr}{\mathbf{r}}
\newcommand{\cc}{\mathbf{c}}
\newcommand{\xx}{\mathbf{x}}
\newcommand{\uu}{\mathbf{u}}
\newcommand{\vv}{\mathbf{v}}
\newcommand{\zz}{\mathbf{z}}
\newcommand{\zero}{\mathbf{0}}
\newcommand{\pp}{\mathbf{p}}
\newcommand{\cG}{\mathcal{G}}
\newcommand{\cT}{\mathcal{T}}
\newcommand{\cQ}{\mathcal{Q}}
\newcommand{\U}{\mathfrak{U}}
\newcommand{\eps}{\varepsilon}
\newcommand{\pos}{\operatorname{pos}}

\title{Additive cubes are avoidable over $\{0,1,2,3\}$, $\{0,1,4\}$ and $\{0,2,5\}$} 

\author{
	\textbf{Damien Jamet}\\
	Université de Lorraine, LORIA\\
	Campus Scientifique, BP 239\\
	F-54506 Vandœuvre-lès-Nancy Cedex, France\\
	\texttt{damien.jamet@loria.fr}
	\and
	\textbf{Matthieu Rosenfeld}\\
	Laboratoire GAATI, Université de la Polynésie française\\
	BP 6570\\
	98702 Faa'a, Polynésie française\\
	\texttt{matthieu.rosenfeld@lirmm.fr}
}

\date{\today}

\begin{document}
	\sloppy
	\maketitle
	
	\begin{abstract}
		We prove that additive cubes are avoidable over the alphabet $\{0,1,2,3\}$, thereby solving a long-standing question initiated by Pirillo and Varricchio in 1994. By extending the reduction and spectral machinery of Cassaigne, Currie, Schaeffer, and Shallit (2014), we prove that the search for additive cubes in such morphic images reduces to a finite graph traversal. We exhibit eight explicit morphisms from $\{0,1,3,4\}^*$ to $\{0,1,2,3\}^*$ of length 3 whose images on the Cassaigne \textit{et al.} infinite word contain no additive cubes.
	\end{abstract}

	\section{Introduction and main result}\label{sec:intro}
	
	Throughout the present document, terminology and notations follow standard conventions on combinatorics on words \cite{Lothaire}. All the alphabets are finite sets included in $\N$. Let $\Sigma \subset \N$ be an alphabet. A finite word $u \in \Sigma^*$ has \emph{length} $|u| = n$. The \emph{sum} of the letters of $u$ is $\som(u) = \sum_{k=0}^{n-1} u_k$ (with $\som(\eps) = 0$).
	
	\begin{definition}[Additive cube]\label{def:cube}
		A non-empty finite word $u$ over an alphabet $\Sigma \subseteq \N$ is an \emph{additive cube} if $u=u_1 u_2 u_3$ with $u_1,\ u_2,\ u_3 \in \Sigma^\star$, $|u_1| = |u_2| = |u_3|$ and $\som(u_1)=\som(u_2)=\som(u_3)$.
		A word (finite or infinite) \emph{is without additive cubes} or \emph{avoids additive cubes}, if none of its finite factors is an additive cube.
	\end{definition}
	
	The question of the existence of an infinite word avoiding additive cubes over a finite alphabet was raised by Pirillo and Varricchio \cite{PV} in 1994. Cassaigne, Currie, Schaeffer and Shallit \cite{CCSS} showed the existence of a pure morphic infinite word over $\{0,1,3,4\}$ without additive cubes. Lietard and Rosenfeld \cite{LR} proved that additive cubes are avoidable by pure morphic words over all alphabets of four integers that are not affinely equivalent to $\{0,1,2,3\}$. 	Rao \cite{Rao} extended the avoidability of additive cubes to ternary alphabets such as $\{0,1,5\}$.
	The question for $\{0,1,2,3\}$ (question~5 of \cite{Rao}) is, to our knowledge, still open.
	
	We exhibit eight morphisms from $\{0,1,3,4\}$ to $\{0,1,2,3\}$ such that the image of the infinite word introduced in~\cite{CCSS} under each of these morphisms avoids additive cubes. 
	
	Our proof scheme is largely inspired by---and builds directly on---the approach of \cite{CCSS}. The key to our method lies in a slightly stronger use of Theorem~\ref{thm:contraction} (Theorem 8 in \cite{CCSS}); although it is restated and slightly reformulated here, this essential property is entirely due to \cite{CCSS}.
	
	From now on, $\Sigma$ denotes the alphabet $\{0,1,3,4\}$, $\varphi: \Sigma^* \to \Sigma^*$  denotes the morphism 
	\[
	\varphi(0)=03,\qquad \varphi(1)=43,\qquad \varphi(3)=1,\qquad \varphi(4)=01 .
	\]
	and $\w = \lim_{n \to \infty} \varphi^n(0) = 0314301103434303101101\ldots = \varphi(\w)$.
	
	\begin{theoreme}[\cite{CCSS}]\label{thm:ccss}
		The infinite word $\w$ contains no additive cube.
	\end{theoreme}
	
	Let $\tau: \Sigma^* \to \{0,1,2,3\}^*$ be a \emph{non-erasing} morphism, that is $|\tau(a)| \ge 1$ for every $a \in \Sigma$. We denote:
	\[
	\bl:= \begin{pmatrix} |\tau(0)| \\ |\tau(1)| \\ |\tau(3)| \\ |\tau(4)| \end{pmatrix} =
	\begin{pmatrix} \ell(0) \\ \ell(1)  \\ \ell(3)  \\ \ell(4) \end{pmatrix}
	\in (\N^*)^4, \quad
	\bs:= \begin{pmatrix} \som(\tau(0)) \\ \som(\tau(1)) \\ \som(\tau(3)) \\ \som(\tau(4)) \end{pmatrix} =
	\begin{pmatrix} s(0) \\ s(1) \\ s(3) \\ s(4) \end{pmatrix}
	\in \N^4, \quad
	\Gamma:= \begin{pmatrix} \transpose{\bl} \\ \transpose{\bs} \end{pmatrix} \in \Z^{2 \times 4},
	\]
	Let $\xx = \tau(\w) = \tau(w_0)\tau(w_1)\tau(w_2)\cdots$. 
	
	We will check that the morphisms $\tau_1, \dots, \tau_8$ listed in Table~\ref{tab:resultats} on page~\pageref{tab:resultats} satisfy two hypotheses:
	\begin{itemize}
		\item \textbf{(H0)} for every $a \in \Sigma$, $\tau(a)$ contains no additive cube;
		\item \textbf{(T)} \emph{transversality}: a specific $4\times 4$ matrix $\Phi$ depending on $\tau$ is invertible (Definition~\ref{def:T}).
	\end{itemize}
	
	\begin{theoreme}[Reduction theorem, precise statement in Theorem~\ref{thm:reduction}]\label{thm:A}
		Let $\tau: \Sigma^* \to \{0,1,2,3\}^*$ be a non-erasing morphism  satisfying (H0) and (T). There exist a finite graph $\cG'$, a starting set $A_0$ of at most $9$ vertices and an accepting set $Z_\tau$ 
		such that, denoting by $R$ the set of vertices of $\cG'$ reachable from $A_0$:
		\[
		\xx = \tau(\w) \text{ avoids additive cubes} \iff R \cap Z_\tau = \emptyset.
		\]
	\end{theoreme}
	
	\begin{theoreme}[Result of the computation, Section~\ref{sec:calcul}]\label{thm:B}
		For each of the eight morphisms $\tau_1,\dots,\tau_8$ of Table~\ref{tab:resultats} on page~\pageref{tab:resultats}, the hypotheses (H0) and (T) are satisfied and the traversal of $\cG'$ gives $R \cap Z_\tau = \emptyset$. By Theorem~\ref{thm:A}, $\tau_i(\w)$ is an infinite word over $\{0,1,2,3\}$ without additive cubes.
	\end{theoreme}

    A direct consequence of Theorem \ref{thm:B} and \cite{Rao} (see Table 3 and Question 5) is  the avoidability of additive cubes over two ternary alphabets:
\begin{corollaire}
	Additive cubes are avoidable over the alphabets $\{0, 1, 4\}$ and $\{0, 2, 5\}$. 
\end{corollaire}

	In the following, we are not going to detail every example, but only focus on a particular case as a running example: $\tau_1: 0 \mapsto 21,\ 1 \mapsto 100,\ 3 \mapsto 3,\ 4 \mapsto 23$. For this morphism, we have $\bl = \transpose{(2,3,1,2)}$, $\bs = \transpose{(3,1,3,5)}$,
	\begin{equation*}
		\Gamma = \begin{pmatrix} 2&3&1&2\\3&1&3&5\end{pmatrix} \text{ and }
		\xx = \tau(0)\,\tau(3)\,\tau(1)\,\tau(4)\,\tau(3)\,\tau(0)\cdots = 21\;3\;100\;23\;3\;21\; \cdots
	\end{equation*}	
	The other morphisms are handled in exactly the same way, and we will provide the numerical results for each of them.
	
	Our entire work is a direct adaptation of the work by Cassaigne, Currie, Schaeffer, and Shallit~\cite{CCSS}. The core proof strategy---analyzing the genealogy of factors, bounding spectral projections, and reducing the search to an induced finite graph---is already present in~\cite{CCSS}. In fact, their original theorem can be viewed as the canonical base case where $\tau = \mathrm{id}_{\Sigma}$, without boundary shifts. Our contribution adapts their machinery to handle the deformations induced by a coding morphism $\tau$.

    The computer verification relying on Theorem~\ref{thm:B} was carried out by a \texttt{C++} program. This program is an updated  version of the code that the first author originally developed to search for candidate morphisms whose fixed points avoid additive cubes. These morphisms were then studied and proven in~\cite{LR}.
    
	\section{Properties of the word $\mathbf{w}$}
	
	While the content of this section is largely based on~\cite{CCSS}, we recall and detail its main components in order to make the present paper self-contained and easy to follow.
	
	\subsection{Notation}
	
	We endow $\R^4$ with the canonical basis $(\ee_0, \ee_1, \ee_3, \ee_4)$ indexed by the letters of $\Sigma$. The \emph{Parikh vector} of a finite word $u$ is $\psi(u) = \transpose{(|u|_0, |u|_1, |u|_3, |u|_4)} \in \N^4$. For $p \le q$, we denote $\w[p,q) = w_p \cdots w_{q-1}$ (with $\w[p,p) = \eps$) and $\sigma(p) = \psi(\w[0,p))$ is the Parikh vector of the prefix $w_0 \cdots w_{p-1}$; hence $\sigma(0)=\zero$ and, for every factor $b = \w[p,q)$, $\psi(b) = \sigma(q)-\sigma(p)$.
	
	The \emph{incidence matrix} of $\varphi$ is defined as the $4 \times 4$ matrix whose columns are the Parikh vectors of the images of $0$, $1$, $3$, and $4$, in that order:
	\[
	M = \left( \psi(\varphi(0)) \;\middle|\; \psi(\varphi(1)) \;\middle|\; \psi(\varphi(3)) \;\middle|\; \psi(\varphi(4)) \right) = \begin{pmatrix}
		1 & 0 & 0 & 1 \\
		0 & 0 & 1 & 1 \\
		1 & 1 & 0 & 0 \\
		0 & 1 & 0 & 0
	\end{pmatrix},
	\qquad \psi(\varphi(u)) = M\psi(u).
	\]
	Let $\eta(p) = |\varphi(\w[0,p))|$ be the length of the image of the prefix $\w[0,p)$ under $\varphi$.
	
	For any integer $q$, the \emph{parent} of $q$, denoted $\parent(q)$, is the unique $p$ such that $\eta(p) \le q < \eta(p+1)$; we then say that $q$ is a \emph{child} of $p$. Roughly speaking, the parent of an index $q$ is the unique index $p$ such that the image $\varphi(w_p)$ covers the position $q$ within the word $\mathbf{w}$.
	
	The word $\w$ begins with $\w[0,20) = 03143011034343031011$ (see Figure~\ref{fig:blocs-parents}). Table~\ref{tab:filiation-valeurs} gives the values of $\eta$, parents, and children for the first few positions.
	
	\begin{figure}[!t]
		\centering
		\begin{tikzpicture}[
			xscale=0.54,
			yscale=1.0,
			lbl/.style={font=\footnotesize\bfseries, text=black!70, anchor=east},
			idx/.style={font=\tiny, text=gray!75!black},
			val/.style={font=\normalsize\bfseries},
			par/.style={font=\footnotesize\bfseries}
			]
			\node[lbl] at (-0.6,  0.48) {indices $q$:};
			\node[lbl] at (-0.6,  0.00) {letters $w_q$:};
			\node[lbl] at (-0.6, -0.58) {parent $t$:};
			\foreach \q/\w in {
				0/0, 1/3, 2/1, 3/4, 4/3, 5/0, 6/1, 7/1, 8/0, 9/3,
				10/4, 11/3, 12/4, 13/3, 14/0, 15/3, 16/1, 17/0, 18/1, 19/1
			} {
				\node[idx] at (\q, 0.48) {\q};
				\node[val] at (\q, 0.00) {\w};
			}
			\foreach \t/\qA/\qB in {
				0/0/1,  1/2/2,   2/3/4,   3/5/6,
				4/7/7,  5/8/9,   6/10/11, 7/12/13,
				8/14/15, 9/16/16, 10/17/18, 11/19/19
			} {
				\draw[semithick, gray!60] 
				(\qA - 0.32, -0.20) -- ++(0, -0.10) 
				-- (\qB + 0.32, -0.30) -- ++(0, 0.10);
				\path (\qA, 0) -- (\qB, 0) coordinate[midway] (mid);
				\node[par] at (mid |- 0, -0.58) {\t};
			}
		\end{tikzpicture}
		\caption{Splitting of the first $20$ positions of $\w$ into images $\varphi(w_t)$. The indices under the brackets indicate the parent $t = \parent(q)$.}
		\label{fig:blocs-parents}
	\end{figure}
	
	\begin{table}[!h]
		\centering\footnotesize
		\caption{Values of $\eta(p)$, $\parent(p)$, and children for the first positions of $\w$.}
		\label{tab:filiation-valeurs}
		\setlength{\tabcolsep}{4.2pt}
		\begin{tabular}{l|cccccccccccc}
			\toprule
			position $p$ & 0 & 1 & 2 & 3 & 4 & 5 & 6 & 7 & 8 & 9 & 10 & 11 \\
			\midrule
			letter $w_p$ & 0 & 3 & 1 & 4 & 3 & 0 & 1 & 1 & 0 & 3 & 4 & 3 \\
			$\eta(p)$ & 0 & 2 & 3 & 5 & 7 & 8 & 10 & 12 & 14 & 16 & 17 & 19 \\
			$\parent(p)$ & 0 & 0 & 1 & 2 & 2 & 3 & 3 & 4 & 5 & 5 & 6 & 6 \\
			\midrule
			children of $p$ & $\{0, 1\}$ & $\{2\}$ & $\{3, 4\}$ & $\{5, 6\}$ & $\{7\}$ & $\{8, 9\}$ & $\{10, 11\}$ & $\{12, 13\}$ & $\{14, 15\}$ & $\{16\}$ & $\{17, 18\}$ & $\{19\}$ \\
			\bottomrule
		\end{tabular}
	\end{table}
	
	In Lemma \ref{lem:eta},	we recall the principal properties of $\eta$ and $\parent$:
	\begin{lemme}[see Section 2.2 in \cite{CCSS}]\label{lem:eta}
		\begin{enumerate}
			\item For every $p$, $\w[\eta(p), \eta(p+1)) = \varphi(w_p)$, that is the children of $p$ are exactly the integers $\eta(p)$, $\eta(p) + 1$, \dots,  $\eta(p) +  |\varphi(w_p)| - 1$.
			\item $\eta$ is strictly increasing, $\eta(0)=0$, and $\eta(p) > p$ for every $p \ge 1$.
			\item $\parent$ is increasing (not necessarily strictly), $\parent(0)=0$, and $\parent(p) < p$ for every $p \ge 1$. In particular, for every $p$, we can find $k \ge 0$ such that $\parent^k(p) = 0$.
		\end{enumerate}
	\end{lemme}
	
	We denote by $\mathcal{A}$ the set of proper prefixes of the images under $\varphi$:
	\[
	\mathcal{A} = \{\eps,\, 0,\, 4\}.
	\]
	
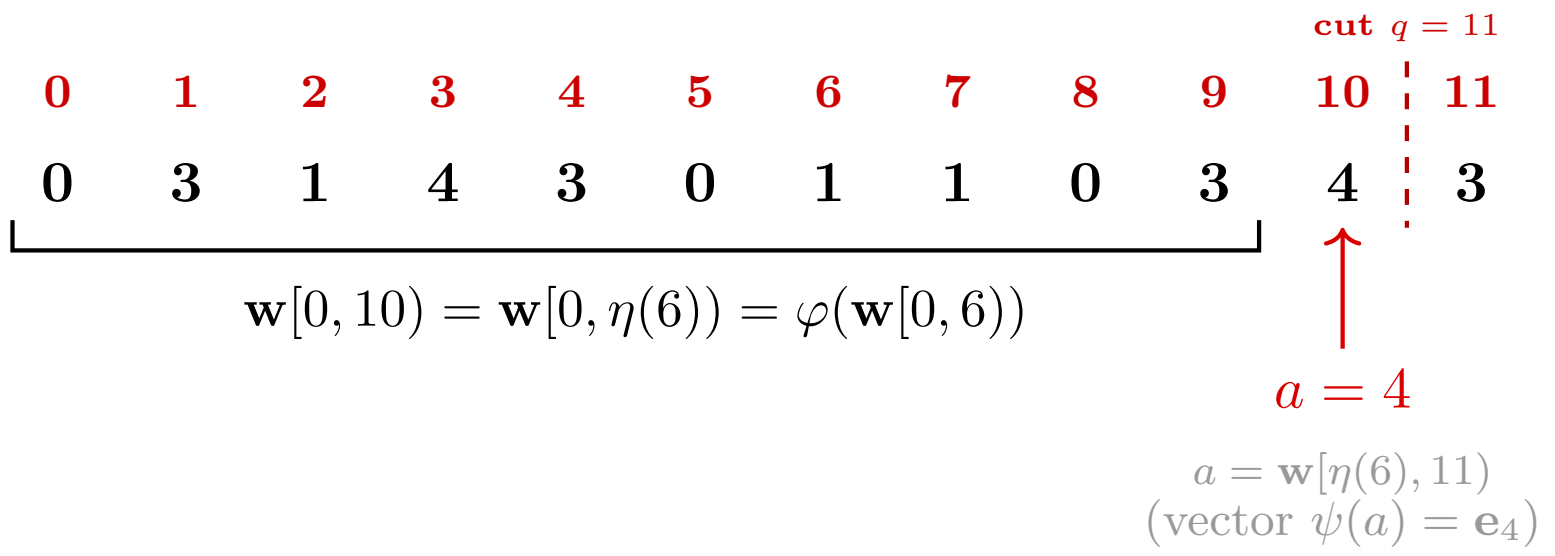
\begin{figure}[!h]
		\centering
		\begin{tikzpicture}[
			xscale=0.85,
			yscale=1.0,
			font=\small,
			idx/.style={font=\footnotesize\bfseries, text=red!80!black},
			val/.style={font=\normalsize\bfseries},
			lbl/.style={font=\small, text=black}
			]
			\foreach \q in {0, 1, 2, 3, 4, 5, 6, 7, 8, 9, 10, 11} {
				\node[idx] at (\q, 0.6) {\q};
			}
			\foreach \q/\w in {
				0/0, 1/3, 2/1, 3/4, 4/3, 5/0, 6/1, 7/1, 8/0, 9/3, 10/4, 11/3
			} {
				\node[val] at (\q, 0.0) {\w};
			}
			\draw[thick] (-0.35, -0.25) -- (-0.35, -0.45) 
			-- (9.35, -0.45) -- (9.35, -0.25);
			\node[lbl, anchor=north] at (4.5, -0.55) {
				$\w[0, 10) = \w[0, \eta(6)) = \varphi(\w[0, 6))$ 
			};
			\draw[<-, thick, red!85!black] (10, -0.30) -- (10, -1.1) 
			node[below, font=\normalsize\bfseries, text=red!85!black] {$a = 4$};
			\node[font=\scriptsize, text=gray!80, anchor=north, align=center] at (10, -1.65) {
				$a = \w[\eta(6), 11)$ \\
				{\footnotesize(vector $\psi(a) = \ee_4$)}
			};
			\draw[dashed, red!70!black, thick] (10.5, 0.8) -- (10.5, -0.3);
			\node[font=\tiny\bfseries, text=red!80!black, anchor=south] at (10.5, 0.8) {cut $q = 11$};
		\end{tikzpicture}
		\caption{Decomposition of the prefix $\w[0, 11)$: the first ten letters form the complete image $\varphi(\w[0,6))$, and the letter $w_{10} = 4$ is a proper prefix $a = \varphi(w_6)[0,1)$ of $\varphi(w_6)$.}
		\label{fig:exemple-lemme-fond}
	\end{figure}

	For every position $q$ with parent $p$, the prefix $\w[0,q)$ is the concatenation of the \textbf{complete image} $\varphi(\w[0,p))$ of the parent prefix and a proper prefix $a \in \mathcal{A}$ of the image of the next letter $w_p$ (Figure~\ref{fig:exemple-lemme-fond}, case $q=11$):
	\[
	\w[0,q) = \varphi(\w[0,p)) \cdot a = \w[0, \eta(p)) \cdot \w[\eta(p), q), \qquad a = \varphi(w_p)[0,i),\quad i = q - \eta(p).
	\]
	
	\begin{lemme}[see Lemma 3 in \cite{CCSS}]\label{lem:fond}
		Let $q \in \N$, $p = \parent(q)$ and $i = q - \eta(p)$ $(0 \le i < |\varphi(w_p)|)$. The word $a = \varphi(w_p)[0, i)$ is a proper prefix of $\varphi(w_p)$ (hence $a \in \mathcal{A}$) and
		\[
		\sigma(q) = M \sigma(p) + \psi(a).
		\]
	\end{lemme}
	
	\begin{corollaire}[see Corollary 4  in \cite{CCSS}]\label{cor:chaine}
		Let $p_0 \in \N$ and $p_{i+1} = \parent(p_i)$. For $i \ge 0$, we set $a_i = \w[\eta(p_{i+1}), p_i) \in \mathcal{A}$. Then
		\[
		\sigma(p_0) = \sum_{i \ge 0} M^i \psi(a_i),
		\]
		a finite sum since $p_i = 0$ and $a_i = \eps$ for $i$ large enough.
	\end{corollaire}
	
	\begin{corollaire}[see Corollary 4  in \cite{CCSS}]\label{cor:facteurs}
		Let $b = \w[p,q)$ be a factor and $b' = \w[\parent(p), \parent(q))$ its \textbf{parent factor}. Then
		\[
		\psi(b) = M \psi(b') + \delta, \qquad \delta = \psi(a') - \psi(a) \in D:= \psi(\mathcal{A}) - \psi(\mathcal{A}),
		\]
		where $a = \w[\eta(\parent(p)), p)$ and $a' = \w[\eta(\parent(q)), q)$ are the proper prefixes of Lemma~\ref{lem:fond} at the boundaries $p$ and $q$.
	\end{corollaire}
	
	\begin{remarque}[Finiteness of boundary perturbations]\label{rem:finitude-D}
		Since $\psi(\eps) = \zero$, $\psi(0) = \ee_0$, $\psi(4) = \ee_4$, the set $D$ contains only $7$ vectors: $D = \{\zero, \pm \ee_0, \pm \ee_4, \pm(\ee_0 - \ee_4)\}$. It is this finiteness that allows us to bound the boundary errors at all scales.
	\end{remarque}
	
	Ancestral chains (parents, grandparents, etc.) are strictly determined by local rules: the edge label (the prefix $a$) depends entirely on the parent letter and the child's index in the image, which allows us to encode this local information in a graph.
	
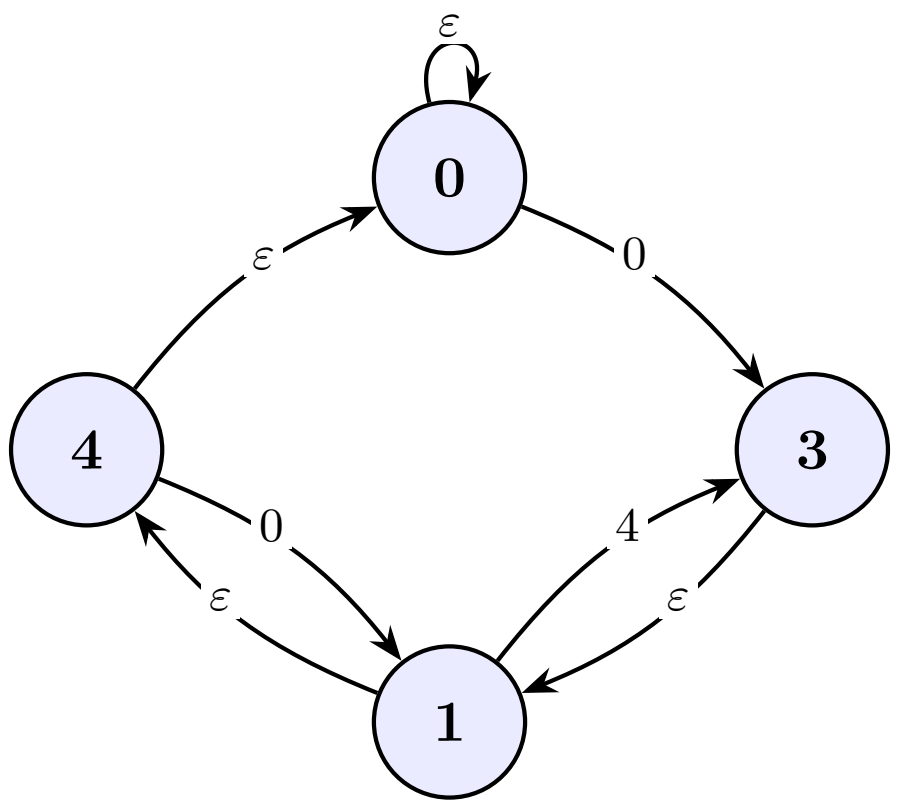
\begin{figure}[t]
		\centering
		\begin{tikzpicture}[
			>=Stealth,
			node distance=2.6cm,
			every node/.style={font=\small},
			state/.style={circle, draw=black, thick, fill=blue!8, minimum size=10mm, font=\normalsize\bfseries},
			lbl/.style={fill=white, inner sep=1.5pt, font=\footnotesize}
			]
			\node[state] (0) at (0, 1.8) {0};
			\node[state] (3) at (2.4, 0) {3};
			\node[state] (1) at (0, -1.8) {1};
			\node[state] (4) at (-2.4, 0) {4};
			\draw[->, thick] (0) edge[loop above, looseness=5] node[lbl, above] {$\eps$} (0);
			\draw[->, thick] (0) to[bend left=15] node[lbl, above left] {$0$} (3);
			\draw[->, thick] (1) to[bend left=15] node[lbl, above left] {$\eps$} (4);
			\draw[->, thick] (1) to[bend left=15] node[lbl, above right] {$4$} (3);
			\draw[->, thick] (3) to[bend left=15] node[lbl, above right] {$\eps$} (1);
			\draw[->, thick] (4) to[bend left=15] node[lbl, above right] {$\eps$} (0);
			\draw[->, thick] (4) to[bend left=15] node[lbl, above left] {$0$} (1);
		\end{tikzpicture}
		\caption{The graph $\cQ$: parent $\to$ child transitions and their labels.}
		\label{fig:graphe-Q}
	\end{figure}

	\begin{definition}[The graph $\cQ$ (Figure \ref{fig:graphe-Q}), See Section 2.4 in \cite {CCSS}]\label{def:Q}	The graph $\cQ$ has the letters of $\Sigma$ as vertices. For each letter $c$ and each $0 \le i < |\varphi(c)|$, it has an edge $e = (c, i)$ with \emph{source} $c$, \emph{target} the letter $\varphi(c)[i]$, and \emph{label} $\varphi(c)[0,i) \in \mathcal{A}$.
	\end{definition}
	
	By Lemma~\ref{lem:eta}(1), if $q$ is a child of $p$, then $q = \eta(p) + i$ with $i < |\varphi(w_p)|$; the pair $(w_p, i)$ is an edge of $\cQ$ with source $w_p$, target $w_q$ and label $a = \w[\eta(p), q)$. In other words, \emph{the children of $p$ in $\w$ are in bijection with the edges of $\cQ$ whose source is $w_p$}.

	\subsection{The spectrum of $M$}\label{sec:spectre}
	
	The characteristic polynomial of $M$ is $X^4 - X^3 - 2X^2 + 2X - 1$ and its roots are:
	\begin{align*}
		\lambda_1 &\approx +1.690284, \\
		\lambda_2 &\approx -1.505068, \\
		\lambda_3 &\approx \phantom{+}0.407392 + 0.476565\,i, \\
		\lambda_4 = \overline{\lambda_3} &\approx \phantom{+}0.407392 - 0.476565\,i \quad (|\lambda_3| = |\lambda_4| \approx 0.626963).
	\end{align*}
	The eigenvalues are all simple, hence $M$ is diagonalizable over $\C$: $M = Q \Lambda Q^{-1}$ with $\Lambda = \operatorname{diag}(\lambda_1, \lambda_2, \lambda_3, \lambda_4)$.
	
	\begin{definition}[Linear forms $\theta_j$]\label{def:theta}
		For each $j$, let $q_j \in \C^4$ be a unit eigenvector ($\|q_j\|_2 = 1$) related to the  eigenvalue $\lambda_j$, normalized such that its first component is real and strictly positive. The \emph{eigencoordinates} $\theta_1, \dots, \theta_4: \C^4 \to \C$ are the rows of $Q^{-1}$:
		\[
		Q = \begin{pmatrix} \vert & \vert & \vert & \vert \\ q_1 & q_2 & q_3 & q_4 \\ \vert & \vert & \vert & \vert \end{pmatrix},
		\qquad
		Q^{-1} = \begin{pmatrix} \text{--- } \theta_1 \text{ ---} \\ \text{--- } \theta_2 \text{ ---} \\ \text{--- } \theta_3 \text{ ---} \\ \text{--- } \theta_4 \text{ ---} \end{pmatrix}.
		\]
		For every $z \in \R^4$: $\theta_j(M z) = \lambda_j \theta_j(z)$, $\theta_4(z) = \overline{\theta_3(z)}$ (hence $|\theta_3(z)| = |\theta_4(z)|$), and $z = \sum_{j=1}^4 \theta_j(z)\, q_j$.
	\end{definition}
	Numerical values:
	\begin{align*}
		Q &\approx \begin{pmatrix}
			0.4724 & 0.1781 & 0.6270 & 0.6270 \\[2pt]
			0.5512 & 0.6714 & -0.2938 - 0.0553\,i & -0.2938 + 0.0553\,i \\[2pt]
			0.6056 & -0.5644 & 0.2782 - 0.4613\,i & 0.2782 + 0.4613\,i \\[2pt]
			0.3261 & -0.4461 & -0.3715 + 0.2988\,i & -0.3715 - 0.2988\,i
		\end{pmatrix}, \\[8pt]
		Q^{-1} &\approx \begin{pmatrix}
			0.5124 & 0.5979 & 0.3537 & 0.6569 \\[3pt]
			0.1806 & 0.6809 & -0.4524 & -0.5724 \\[3pt]
			0.5788 - 0.5749\,i & -0.3219 + 0.2183\,i & -0.0690 + 0.6165\,i & -0.1662 - 0.6810\,i \\[3pt]
			0.5788 + 0.5749\,i & -0.3219 - 0.2183\,i & -0.0690 - 0.6165\,i & -0.1662 + 0.6810\,i
		\end{pmatrix}.
	\end{align*}
	
	\subsection{Contracting directions: every deviation is bounded}
	We now exploit the contracting eigendirections to bound any Parikh deviation between factors of $\w$.
	
	\begin{definition}[Set of paths $\mathcal{D}_\ell$, Definition 7 in \cite{CCSS}]\label{def:D-ell}
		For $\ell \ge 1$, $\mathcal{D}_\ell \subset \Z^4$ is the set of vectors $\sum_{i=0}^{\ell-1} M^i \psi(a_i)$, where $(e_0, \dots, e_{\ell-1})$ ranges over the sequences of edges of $\cQ$ such that the target of $e_{i+1}$ is the source of $e_i$ for $0 \le i < \ell - 1$ (the path is read from the oldest generation $e_{\ell-1}$ to the most recent one $e_0$), and $a_i$ is the label of $e_i$.
	\end{definition}
	
	The graph $\cQ$ has $466$ paths of length $9$, which generate $301$ distinct vectors in $\mathcal{D}_9$ (see proof of Theorem 8 \cite{CCSS}). We set:
	\[
	C_* = \max_{s,t \in \mathcal{A}} |\theta_3(\psi(s) - \psi(t))|, \qquad
	C_3 = 2 \left( \max_{\alpha, \alpha' \in \mathcal{D}_9} |\theta_3(\alpha) - \theta_3(\alpha')| + \frac{C_* \, |\lambda_3|^9}{1 - |\lambda_3|} \right) \approx 2.175815.
	\]
	
	\begin{theoreme}[Theorem 8 in \cite{CCSS}]\label{thm:contraction}
		For every pair of factors $b, c$ of $\w$ (not necessarily of the same length, not necessarily consecutive):
		\[
		|\theta_3(\psi(b) - \psi(c))| \le C_3 \qquad\text{and}\qquad |\theta_4(\psi(b) - \psi(c))| \le C_3.
		\]
	\end{theoreme}
	
	\subsection{Expanding directions: stability under ascent}
	
	Let $b$, $c$ be two factors of $\w$ with parents $b'$, $c'$. By Corollary~\ref{cor:facteurs}, the deviation $\Delta = \psi(b) - \psi(c)$ and the deviation $\Delta' = \psi(b') - \psi(c')$ satisfy:
	\[
	\Delta = M \Delta' + \boldsymbol{\rho}, \qquad \boldsymbol{\rho} = \delta(b) - \delta(c) \in E:= D - D,
	\]
	where $E$ is a finite subset of $\Z^4$. In the expanding directions ($j=1,2$), applying this relation towards the children would amplify the deviation by a factor $|\lambda_j| > 1.5$ at each generation. We invert it: applying $\theta_j$,
	\[
	\theta_j(\Delta') = \frac{\theta_j(\Delta) - \theta_j(\boldsymbol{\rho})}{\lambda_j}.
	\]
	
	\begin{definition}[Constants  $\kappa_1$ and  $\kappa_2$]\label{def:kappa}
		For $j \in \{1,2\}$, $\kappa_j = \max_{\boldsymbol{\rho} \in E} |\theta_j(\boldsymbol{\rho})|$; numerically $\kappa_1 \approx 1.3137$ and $\kappa_2 \approx 1.5060$.
	\end{definition}
	
	Iterating from an initial deviation $\Delta_0$, the deviation at the $i$-th generation is:
	\[
	\Delta_i = M^{-i} \Delta_0 - \sum_{k=0}^{i-1} M^{-(i-k)} \boldsymbol{\rho}_k, \qquad\text{hence}\qquad \theta_j(\Delta_i) = \frac{\theta_j(\Delta_0)}{\lambda_j^i} - \sum_{m=1}^{i} \frac{\theta_j(\boldsymbol{\rho}_{i-m})}{\lambda_j^m}.
	\]
	
	\begin{lemme}[Stability under ascent]\label{lem:remontee}
		Let $j \in \{1,2\}$ and $(\Delta_i)_{i \ge 0}$ be the sequence of deviations along a chain of ancestors, with $\boldsymbol{\rho}_k \in E$. If $|\theta_j(\Delta_0)| \le B$, then for every $i \ge 0$:
		\[
		|\theta_j(\Delta_i)| \le \max\left( B, \; \frac{\kappa_j}{|\lambda_j| - 1} \right).
		\]
	\end{lemme}
	
	\begin{proof}
		The case $i=0$ holds. For $i \ge 1$, the triangle inequality gives
		\[
		|\theta_j(\Delta_i)| \le \frac{B}{|\lambda_j|^i} + \kappa_j \sum_{m=1}^i \frac{1}{|\lambda_j|^m} = \frac{B}{|\lambda_j|^i} + \frac{\kappa_j}{|\lambda_j| - 1} \left( 1 - \frac{1}{|\lambda_j|^i} \right).
		\]
		Since $|\lambda_j| > 1$, $\alpha = |\lambda_j|^{-i} \in\, ]0,1]$ and the right-hand side is a convex combination $\alpha B + (1-\alpha)\frac{\kappa_j}{|\lambda_j|-1}$ of $B$ and of $\frac{\kappa_j}{|\lambda_j|-1}$, hence bounded by their maximum.
	\end{proof}
	
	Consequently, at every scales, the Parikh deviation between any two factors of $\w$ is bounded: by $C_3$ in directions $3,4$ (Theorem~\ref{thm:contraction}), and by $\max(B, \kappa_j/(|\lambda_j|-1))$ in directions $1,2$ (Lemma~\ref{lem:remontee}), where $B$ is the initial bound.
	
	\section{Additive cubes in $\xx = \tau(\w)$ via $\w$}
	
	\subsection{Positions in $\xx$}
	
	For $a \in \Sigma$ and $0 \le r \le \ell(a)$, we denote $h_r(a) = \som(\tau(a)[0,r)) = \sum_{k=0}^{r-1} \tau(a)_k$ the sum of the prefix of length $r$ of $\tau(a)$ (with $h_0(a) = 0$).
	
	In general, a position $q$ of $\xx$ does not fall strictly between two images, but occurs at some rank (offset) $r$ within the image of a letter $w_p$. There exists a unique pair $(p,r)$ with $0 \le r < \ell(w_p)$ and
	\[
	\xx[0,q) = \tau(\w[0,p)) \cdot \tau(w_p)[0,r).
	\]
	The offset $r$ is the number of letters ``already consumed" from the image $\tau(w_p)$.
	
	\begin{lemme}[Positions]\label{lem:pos-tau}
		Let $P = \{(p,r): p \in \N,\ 0 \le r < \ell(w_p)\}$ and $\pos(p,r) = \bl \cdot \sigma(p) + r$. The map $\pos: P \to \N$ is a bijection, strictly increasing for the lexicographic order, and the sum $H(q) = \som(\xx[0,q))$ satisfies, for $q = \pos(p,r)$,
		\[
		H(q) = \bs \cdot \sigma(p) + h_r(w_p).
		\]
	\end{lemme}
	
	\begin{proof}
		The prefix $\tau(\w[0,p))$ of $\xx$ has length $\bl \cdot \sigma(p)$ and sum $\bs \cdot \sigma(p)$; the next $r$ letters form a prefix of $\tau(w_p)$ of length $r < \ell(w_p)$ and sum $h_r(w_p)$. The monotonicity comes from the fact that $\bl \cdot \sigma(p+1) = \bl \cdot \sigma(p) + \ell(w_p)$. Surjectivity comes from the fact that every position $q$ can be written in this way (see above).
	\end{proof}
	
	\subsection{The discrete curve of cumulative sums}
	
	For $q = \pos(p,r)$, we denote
	\begin{equation}\label{eq:decompo-lambda}
		\Lambda(q) = \begin{pmatrix} q \\ H(q) \end{pmatrix}
		= \underbrace{\begin{pmatrix} \bl \cdot \sigma(p) \\ \bs \cdot \sigma(p) \end{pmatrix}}_{= \Gamma\sigma(p):\ \text{whole images}} + \underbrace{\begin{pmatrix} r \\ h_r(w_p) \end{pmatrix}}_{\text{offset within the current image}}.
	\end{equation}
	For $q_0 = \pos(p_0,r_0) \le q_1 = \pos(p_1,r_1)$ and $d_0 = \w[p_0,p_1)$:
	\[
	\Lambda(q_1) - \Lambda(q_0) = \begin{pmatrix} |\xx[q_0,q_1)| \\ \som(\xx[q_0,q_1)) \end{pmatrix}
	= \Gamma\,\psi(d_0) + \begin{pmatrix} r_1 - r_0 \\ h_{r_1}(w_{p_1}) - h_{r_0}(w_{p_0}) \end{pmatrix}.
	\]
	
	\subsection{Equations of a cube}
	
	Let $0 \le q_0 < q_1 < q_2 < q_3$, $q_j = \pos(p_j, r_j)$, $c_j = w_{p_j}$ (hence $p_0 \le p_1 \le p_2 \le p_3$ by Lemma~\ref{lem:pos-tau}), and $d_j = \w[p_j, p_{j+1})$ for $j = 0,1,2$. We set:
	\[
	\uu = \psi(d_1) - \psi(d_0) = \sigma(p_0) - 2\sigma(p_1) + \sigma(p_2), \qquad
	\vv = \psi(d_2) - \psi(d_1) = \sigma(p_1) - 2\sigma(p_2) + \sigma(p_3),
	\]
	the \emph{imbalances} between neighboring blocks of $\w$. The three factors $\xx[q_0,q_1)$, $\xx[q_1,q_2)$, $\xx[q_2,q_3)$ form an additive cube if and only if $\Lambda(q_1) - \Lambda(q_0) = \Lambda(q_2) - \Lambda(q_1) = \Lambda(q_3) - \Lambda(q_2)$, that is, if and only if:
	\[
	\Lambda(q_0) - 2\Lambda(q_1) + \Lambda(q_2) = \zero \quad\text{and}\quad \Lambda(q_1) - 2\Lambda(q_2) + \Lambda(q_3) = \zero.
	\]
	By~\eqref{eq:decompo-lambda}, the first expression equals $\Gamma\uu + \bb(\cc,\rr)$ and the second $\Gamma\vv + \bb'(\cc,\rr)$, with:
	\begin{align*}
		\bb(\cc,\rr) &= \begin{pmatrix} r_0 - 2r_1 + r_2 \\ h_{r_0}(c_0) - 2h_{r_1}(c_1) + h_{r_2}(c_2) \end{pmatrix}, &
		\bb'(\cc,\rr) &= \begin{pmatrix} r_1 - 2r_2 + r_3 \\ h_{r_1}(c_1) - 2h_{r_2}(c_2) + h_{r_3}(c_3) \end{pmatrix}.
	\end{align*}
	Thus:
	\begin{equation}\label{eq:systeme-bords}
		\xx[q_0, q_1) \cdot \xx[q_1, q_2) \cdot \xx[q_2, q_3) \text{ is an additive cube} \iff \Gamma \uu + \bb(\cc,\rr) = \Gamma \vv + \bb'(\cc,\rr) = \zero.
	\end{equation}
	The terms $\bb,\bb'$ depend only on the four letters $\cc = (c_0,c_1,c_2,c_3)$ and the four offsets $\rr = (r_0,r_1,r_2,r_3)$: they take only finitely many values.
	
	\subsection{Hypothesis (H0) and accepting set}
	
	\begin{definition}[Local hypothesis (H0)]\label{def:H0}
		The morphism $\tau$ satisfies (H0) if for every $a \in \Sigma$, $\tau(a)$ contains no additive cube;
	\end{definition}
	
	The morphisms in Table~\ref{tab:resultats} on page~\pageref{tab:resultats} have images of length at most $3$ and no image is of the form $ccc$. Hence every morphism satisfies hypothesis (H0).
	
	\begin{definition}[Accepting set $Z_\tau$]\label{def:Ztau}
		The set $Z_\tau \subset \Sigma^4 \times \Z^4 \times \Z^4$ is the set of triples $(\cc,\uu,\vv)$ for which there exist offsets $r_j \in \{0, \dots, \ell(c_j) - 1\}$ ($j = 0,\dots,3$) such that
		\[
		\Gamma\uu + \bb(\cc,\rr) = \zero \quad\text{and}\quad \Gamma\vv + \bb'(\cc,\rr) = \zero.
		\]
	\end{definition}
	
	Then, an additive cube of $\xx$ with boundaries $(p_j,r_j)$ corresponds to a triple $(\cc,\uu,\vv) \in Z_\tau$, where $\cc = (w_{p_0}, w_{p_1}, w_{p_2}, w_{p_3})$ is the quadruple of letters at the boundaries, while $\uu$ and $\vv$ measure the Parikh imbalances between the consecutive blocks $d_j = \w[p_j, p_{j+1})$. The offsets $r_j$ are only required    in this final membership test.
	
	\section{The tree of quadruples and the descent automaton}\label{sec:automate}
	
	\subsection{The tree $\cT$ of quadruples, the map $g$ and the automaton $\cG$}
	
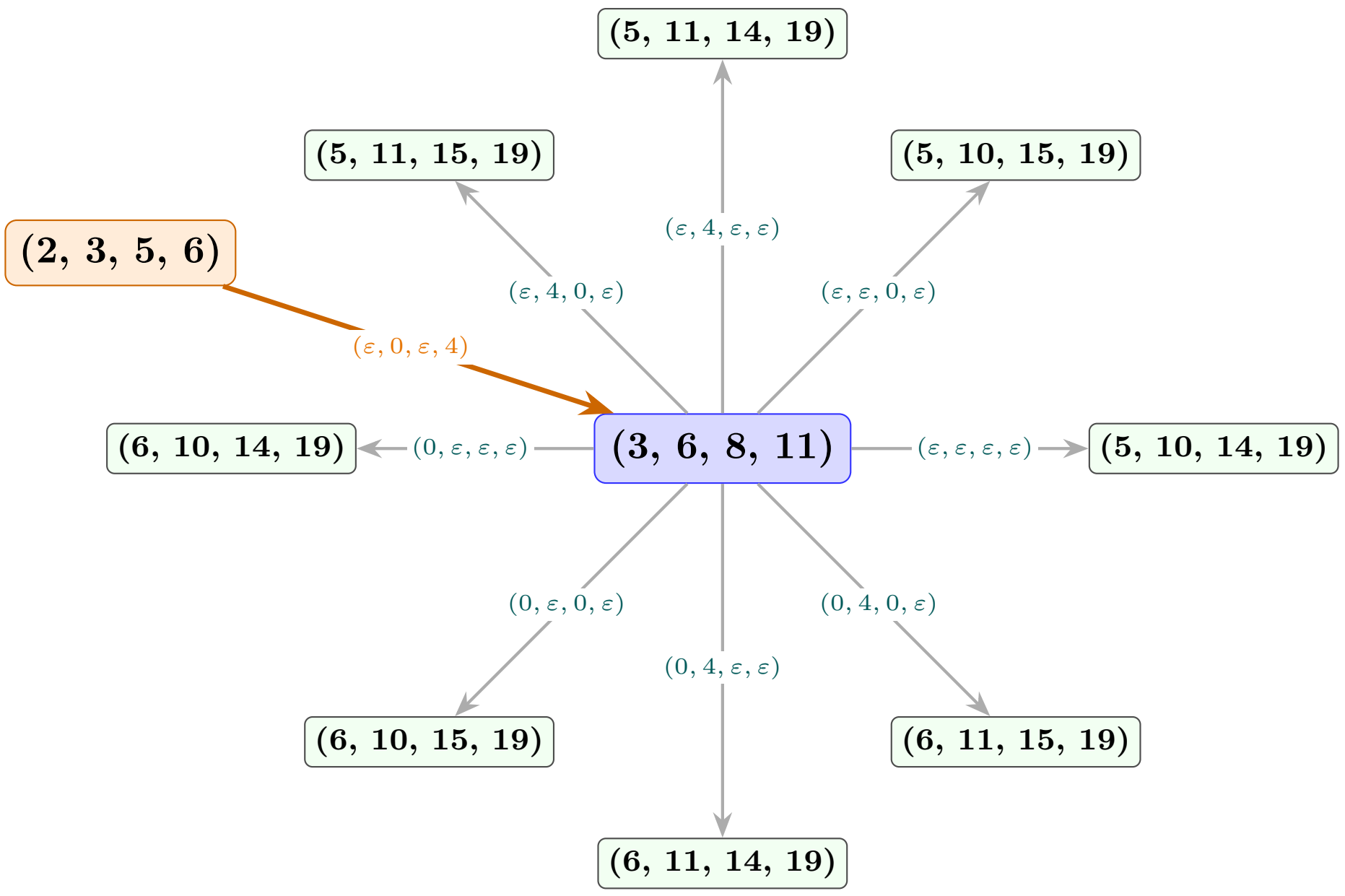
\begin{figure}[t]
		\centering
		\begin{tikzpicture}[
			>=Stealth,
			root/.style={rectangle, draw=blue!80, fill=blue!15, rounded corners=3pt, inner sep=4pt, font=\small\bfseries},
			parent/.style={rectangle, draw=orange!80!black, fill=orange!15, rounded corners=3pt, inner sep=3.5pt, font=\small\bfseries},
			enfant/.style={rectangle, draw=black!70, fill=green!5, rounded corners=2pt, inner sep=2.5pt, font=\scriptsize\bfseries},
			elbl/.style={fill=white, inner sep=1.2pt, font=\tiny\bfseries, text=teal!70!black},
			elblP/.style={fill=white, inner sep=1.5pt, font=\tiny\bfseries, text=orange!90!black}
			]
			\node[root] (root) at (0,0) {(3, 6, 8, 11)};
			\node[parent] (par) at (162:5.8cm) {(2, 3, 5, 6)};
			\draw[->, line width=1.3pt, orange!80!black] (par) -- node[elblP, pos=0.48] {$(\eps, 0, \eps, 4)$} (root);
			\foreach \i/\lbl/\a/\r in {
				0/{(5, 10, 14, 19)}/{(\eps, \eps, \eps, \eps)}/{4.5cm},
				45/{(5, 10, 15, 19)}/{(\eps, \eps, 0, \eps)}/{3.8cm},
				90/{(5, 11, 14, 19)}/{(\eps, 4, \eps, \eps)}/{3.8cm},
				135/{(5, 11, 15, 19)}/{(\eps, 4, 0, \eps)}/{3.8cm},
				180/{(6, 10, 14, 19)}/{(0, \eps, \eps, \eps)}/{4.5cm},
				225/{(6, 10, 15, 19)}/{(0, \eps, 0, \eps)}/{3.8cm},
				270/{(6, 11, 14, 19)}/{(0, 4, \eps, \eps)}/{3.8cm},
				315/{(6, 11, 15, 19)}/{(0, 4, 0, \eps)}/{3.8cm}
			} {
				\node[enfant] (n\i) at (\i:\r) {\lbl};
				\draw[->, thick, gray!65] (root) -- node[elbl, pos=0.52] {$\a$} (n\i);
			}
		\end{tikzpicture}
		\caption{A fragment of the tree $\cT$: the parent quadruple $(2, 3, 5, 6)$ generates $(3, 6, 8, 11)$ with the quadruple of labels $(\eps, 0, \eps, 4)$. The latter has $2^3 = 8$ children because three of its positions (those equal to $4$, $1$ and $0$ in $\w$) have an image of length $2$, whereas the position equal to $3$ has an image of length $1$.}
		\label{fig:arbre-etoile}
	\end{figure}

	An additive cube is delimited by \emph{four} positions. We therefore consider quadruples $\pp = (p_0,p_1,p_2,p_3) \in \N^4$ of positions. A \emph{child} of $\pp$ is a quadruple $\pp' := (p'_0,p'_1,p'_2,p'_3)$ where each $p'_j$ is a child of $p_j$. The \emph{parent} of $\pp$ is $\parent(\pp) := (\parent(p_0),\dots,\parent(p_3))$. This defines a directed tree $\cT$ rooted at $(0,0,0,0)$ with edges from parents to children (see Figure~\ref{fig:arbre-etoile} for an example). Each quadruple has at most $2^4 = 16$ children. We \emph{do not require} the quadruples $\pp$ to be ordered: relative order between coordinates can be modified during the tree descent (from parents to children) when two coordinates in the parent are equal, but it will be fully recovered upon arrival by the equations~\eqref{eq:systeme-bords}.
	
	For instance, the following edges exist in $\cT$ 
	\begin{align*}
		(3, 3, 5, 6) &\xrightarrow{(\eps,\, 0,\, \eps,\, \eps)} (5, 6, 8, 10) \\ 
		(3, 3, 5, 6) &\xrightarrow{(0,\, \eps,\, \eps,\, \eps)} (6, 5, 8, 10) \\ 
		(3, 3, 5, 6) &\xrightarrow{(0,\, \eps,\, 0,\, 4)} (6, 5, 9, 11) 
	\end{align*}
	
	Let $\pp \in \N^4$ and let $\pp'$ be a child of $\pp$. By Lemma~\ref{lem:fond}, for each $j$, one has $\sigma(p'_j) = M\sigma(p_j) + \psi(a_j)$, where $a_j = \varphi(w_{p_j})[0, p'_j - \eta(p_j)) \in \mathcal{A}$ is the label of the edge taken in $\cQ$. Combining these equalities, the imbalances $\uu' = \sigma(p'_0) - 2\sigma(p'_1) + \sigma(p'_2)$ and $\vv' = \sigma(p'_1) - 2\sigma(p'_2) + \sigma(p'_3)$ are equal to
	\[
	\uu' = M \uu + \psi(a_0) - 2\psi(a_1) + \psi(a_2), \qquad
	\vv' = M \vv + \psi(a_1) - 2\psi(a_2) + \psi(a_3).
	\]
	They therefore depend neither on $\pp$ nor on $\pp'$, but only on $(\uu,\vv)$ and the corresponding labels in $\cQ$. We define the state map $g$ as follows:
	\begin{equation}\label{eq:def_g}
		g(\pp) = \bigl( (w_{p_0}, w_{p_1}, w_{p_2}, w_{p_3}), \; \sigma(p_0) - 2\sigma(p_1) + \sigma(p_2), \; \sigma(p_1) - 2\sigma(p_2) + \sigma(p_3) \bigr) \in \Sigma^4 \times \Z^4 \times \Z^4.
	\end{equation}
	In particular, a quadruple of positions $\pp$ delimits an additive cube in $\xx = \tau(\w)$ (for some offsets $\rr$) if and only if $g(\pp) \in Z_\tau$.
	
	Instead of storing the positions $\pp = (p_0, p_1, p_2, p_3)$, which grow arbitrarily large, we only need to remember the four boundary letters and the pair of imbalance vectors  $(\uu, \vv)$. The automaton $\cG$ (see Definition \ref{def:automate-G} below) computes their evolution during the descent by multiplying by $M$ and adding the boundary corrections. Searching for an additive cube thus reduces to a reachability problem in $\cG$.
	
	\begin{definition}[The infinite automaton $\cG$]\label{def:automate-G}
		The vertices of $\cG$ are the states $(\cc,\uu,\vv) \in \Sigma^4 \times \Z^4 \times \Z^4$. For each letter $c_j$ of the quadruple $\cc$, let $e_j = (c_j, i_j)$ be an edge of $\cQ$ with source $c_j$, target $c'_j$, and label $a_j \in \mathcal{A}$. Then $\cG$ contains a directed edge from $(\cc,\uu,\vv)$ to $(\cc',\uu',\vv')$, where
		\[
		\uu' = M \uu + \psi(a_0) - 2\psi(a_1) + \psi(a_2), \qquad
		\vv' = M \vv + \psi(a_1) - 2\psi(a_2) + \psi(a_3).
		\]
	\end{definition}
	
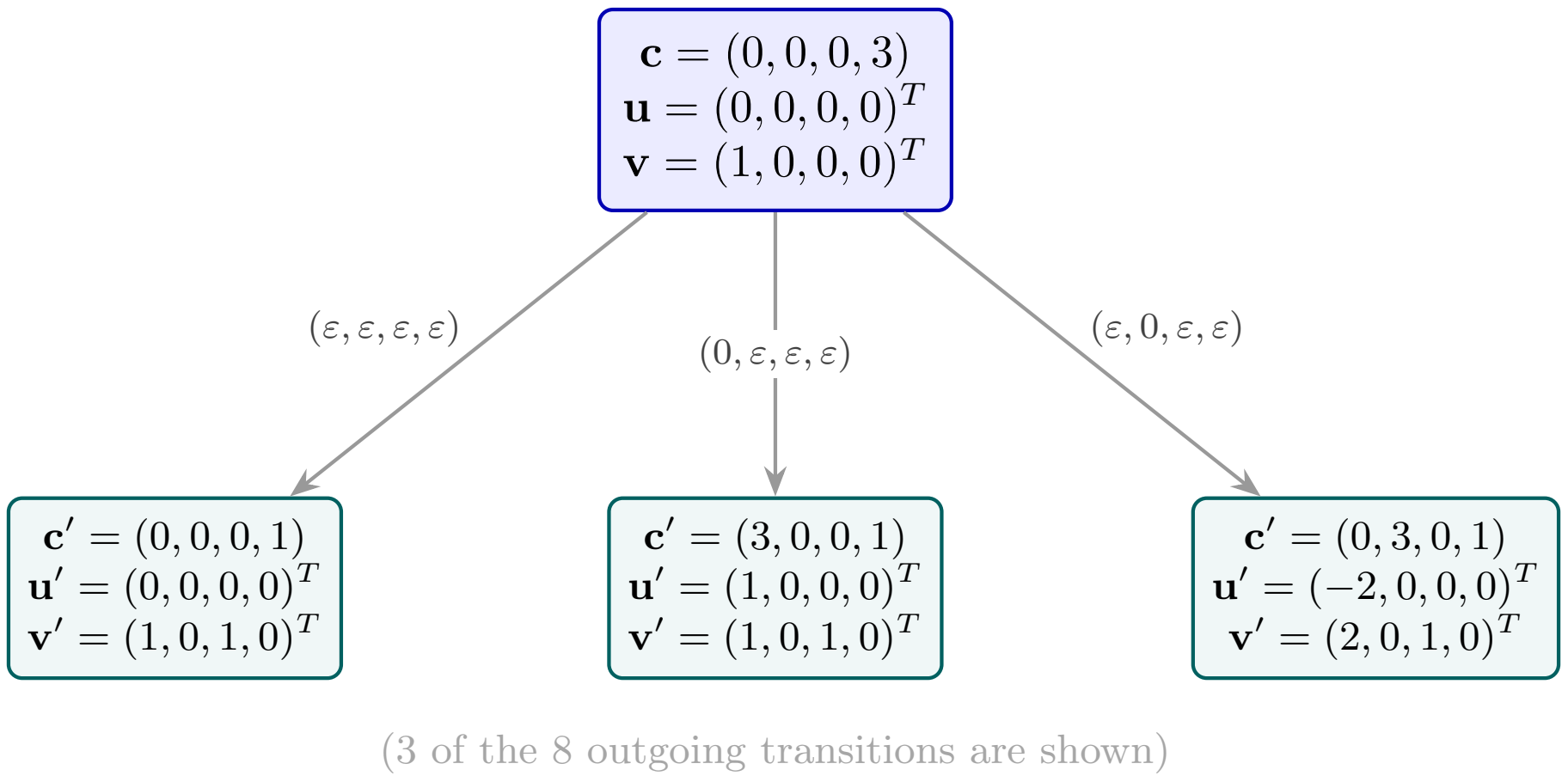
\begin{figure}[!h]
		\centering
		\begin{tikzpicture}[
			>=Stealth,
			xscale=1.05,
			yscale=1.15,
			state/.style={rectangle, draw=blue!70!black, thick, fill=blue!8, rounded corners=3pt, inner sep=5pt, font=\small, align=center},
			childstate/.style={rectangle, draw=teal!75!black, thick, fill=teal!6, rounded corners=3pt, inner sep=4pt, font=\footnotesize, align=center},
			edge lbl/.style={fill=white, inner sep=1.2pt, font=\scriptsize\bfseries, text=black!75}
			]
			\node[state] (root) at (0, 0) {
				$\cc = (0, 0, 0, 3)$ \\
				$\uu = (0, 0, 0, 0)^T$ \\
				$\vv = (1, 0, 0, 0)^T$
			};
			
			\node[childstate] (c1) at (-4.4, -3.2) {
				$\cc' = (0, 0, 0, 1)$ \\
				$\uu' = (0, 0, 0, 0)^T$ \\
				$\vv' = (1, 0, 1, 0)^T$
			};
			\node[childstate] (c2) at (0, -3.2) {
				$\cc' = (3, 0, 0, 1)$ \\
				$\uu' = (1, 0, 0, 0)^T$ \\
				$\vv' = (1, 0, 1, 0)^T$
			};
			\node[childstate] (c3) at (4.4, -3.2) {
				$\cc' = (0, 3, 0, 1)$ \\
				$\uu' = (-2, 0, 0, 0)^T$ \\
				$\vv' = (2, 0, 1, 0)^T$
			};
			
			\draw[->, thick, gray!80] (root) -- node[edge lbl, above left] {$(\eps, \eps, \eps, \eps)$} (c1);
			\draw[->, thick, gray!80] (root) -- node[edge lbl, fill=white] {$(0, \eps, \eps, \eps)$} (c2);
			\draw[->, thick, gray!80] (root) -- node[edge lbl, above right] {$(\eps, 0, \eps, \eps)$} (c3);
			
			\node[font=\footnotesize, text=gray!70] at (0, -4.3) {(3 of the 8 outgoing transitions are shown)};
		\end{tikzpicture}
		\caption{Local transitions in $\cG$ from the state $((0,0,0,3), \zero, \ee_0)$. Each edge corresponds to a 4-tuple of transitions in $\cQ$ with labels $(a_0,a_1,a_2,a_3)$, from which $\uu'$ and $\vv'$ are computed using $M$, $\uu$, $\vv$, and the Parikh vectors $\psi(a_0)$, $\psi(a_1)$, $\psi(a_2)$ and $\psi(a_3)$.\label{fig:transitions-locales-G}}
	\end{figure}

	Hence, a state has as many outgoing transitions as there are combinations of choices for its four letters in $\cQ$ (at most $2^4 = 16$). Figure~\ref{fig:transitions-locales-G} illustrates three transitions from the state $(\cc, \uu, \vv) = ((0,0,0,3), \zero, \ee_0)$. Since the letter $0$ admits two outgoing edges in $\cQ$ (towards $0$ with label $\eps$, or towards $3$ with label $0$) and the letter $3$ has only one (towards $1$ with label $\eps$), this state has exactly $2^3 \times 1 = 8$ outgoing edges in $\cG$.
	
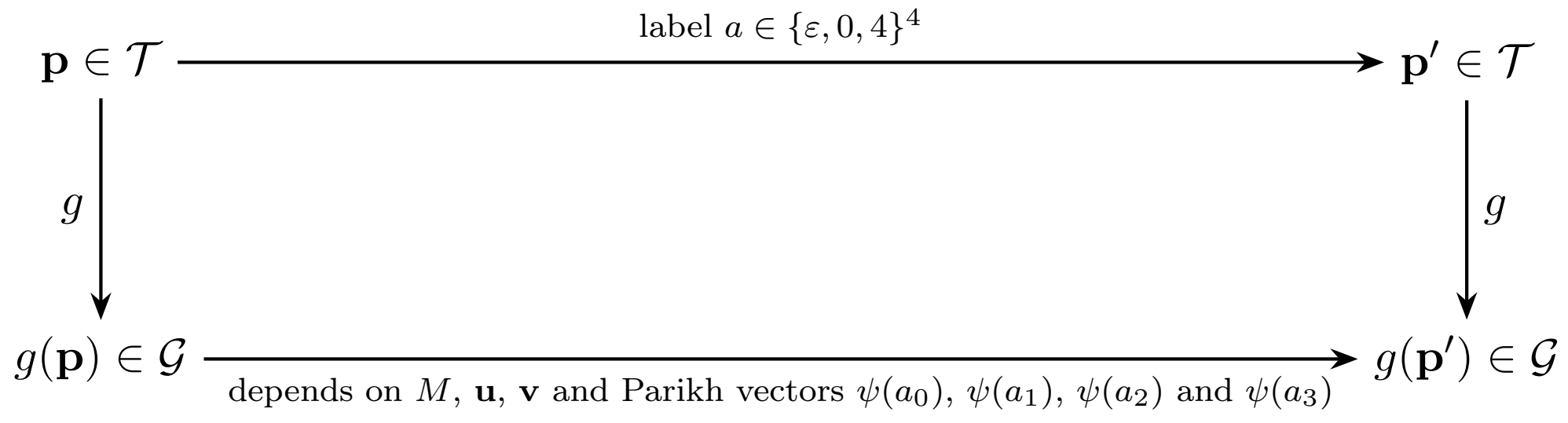
\begin{figure}[!h]
		\centering
		\begin{tikzpicture}[>=Stealth, auto]
			\node (T_p)  at (0, 2.5) {$\pp \in \cT$};
			\node (T_pp) at (11.5, 2.5) {$\pp' \in \cT$};
			\node (G_p)  at (0, 0)   {$g(\pp) \in \cG$};
			\node (G_pp) at (11.5, 0)   {$g(\pp') \in \cG$};
			
			\draw[->, thick] (T_p)  -- node[above, font=\scriptsize] {label $a \in \{\eps, 0, 4\}^4$} (T_pp);
			\draw[->, thick] (G_p)  -- node[below, font=\scriptsize] {depends on $M$, $\uu$, $\vv$ and Parikh vectors $\psi(a_0)$, $\psi(a_1)$, $\psi(a_2)$ and $\psi(a_3)$} (G_pp);
			
			\draw[->, thick] (T_p)  -- node[left]  {$g$} (G_p);
			\draw[->, thick] (T_pp) -- node[right] {$g$} (G_pp);
		\end{tikzpicture}
		\caption{Local commutation between tree descent in $\cT$ and transitions in $\cG$.}
		\label{fig:diagramme-commutatif-G}
	\end{figure}

	\begin{lemme}[Compatibility of $\cG$ with tree descent]\label{lem:automate-marche}
		Let $\pp \in \N^4$.
		\begin{enumerate}
			\item If $\pp'$ is a child of $\pp$ in $\cT$, then $g(\pp) \to g(\pp')$ is an edge in $\cG$.
			\item Conversely, every edge in $\cG$ leaving $g(\pp)$ has the form $g(\pp) \to g(\pp')$ for some child $\pp'$ of $\pp$ in $\cT$.
		\end{enumerate}
		Consequently, the map $g$ establishes a label-preserving bijection between paths from $\pp$ in $\cT$ and paths from $g(\pp)$ in $\cG$ (see Figure~\ref{fig:diagramme-commutatif-G}).
	\end{lemme}
	
	\begin{proof}
		\begin{enumerate}
			\item For each $j \in \{0,1,2,3\}$, let $p'_j = \eta(p_j) + i_j$ with $0 \le i_j < |\varphi(w_{p_j})|$. By Lemma~\ref{lem:eta}(1), this gives the edge $e_j = (w_{p_j}, i_j)$ of $\cQ$ with target $w_{p'_j}$ and label $a_j$. Applying Lemma~\ref{lem:fond} to each coordinate gives the formulas for $\uu'$ and $\vv'$ in Definition~\ref{def:automate-G}, so $g(\pp) \to g(\pp')$ is an edge of $\cG$.
			\item  Conversely, an edge leaving $g(\pp) = (\cc,\uu,\vv)$ corresponds to four edges $e_j = (c_j, i_j)$ in $\cQ$. Setting $p'_j = \eta(p_j) + i_j$, each $p'_j$ is a child of $p_j$ by Lemma~\ref{lem:eta}(1), and part~(1) shows that the target in $\cG$ is precisely $g(\pp')$.
		\end{enumerate}
	\end{proof}
	
	\subsection{The $9$ starting states}\label{sec:Germes}
	
	Let $C$ be an additive cube in $\xx = \tau(\w)$ with boundary positions $q_j = \pos(p_j,r_j)$ satisfying $p_0 \le p_1 \le p_2 \le p_3$ and $p_0 < p_3$ (see Proposition~\ref{prop:reduction}). Following the ancestral chain of $\pp = (p_0, p_1, p_2, p_3)$ in $\w$, there exists $k \ge 0$ such that the four coordinates of $\parent^k(\pp)$ fall within the image $\varphi(w_t)$ of a single letter $w_t$. Since $|\varphi(a)| \le 2$ for every $a \in \Sigma$, these four positions can only take two consecutive values, say $e$ and $e+1$. With $p_0 < p_3$, this forces $|\varphi(w_t)| = 2$ (so $w_t \in \{0, 1, 4\}$) and leaves only three possible non-decreasing configurations for $\parent^k(\pp)$:
	\[
	(e, e, e, e+1), \quad (e, e, e+1, e+1), \quad \text{and} \quad (e, e+1, e+1, e+1).
	\]
	Since the parent letter $w_t$ belongs to $\{0, 1, 4\}$, this produces exactly $3 \times 3 = 9$ candidate starting quadruples.
	
	\begin{definition}[Seed set $X$]\label{def:gX}
		We define the set of seeds $X \subset \N^4$ by
		\[
		X := \left\{ \pp \in \N^4 \;\middle|\; p_0 \le p_1 \le p_2 \le p_3,\; p_0 < p_3,\text{ and } \parent(p_0) =  \parent(p_1) =  \parent(p_2) =  \parent(p_3) \right\}.
		\]
	\end{definition}
	
	\begin{lemme}[The nine starting states]\label{lem:gX}
		The image set $g(X)$ consists of the nine vertices given, for each $xy \in \{03, 43, 01\}$, by:
		\[
		\bigl( (x,x,x,y),\, \zero,\, \psi(x) \bigr), \qquad
		\bigl( (x,x,y,y),\, \psi(x),\, -\psi(x) \bigr), \qquad
		\bigl( (x,y,y,y),\, -\psi(x),\, \zero \bigr).
		\]
	\end{lemme}
	
	\begin{proof}
		Let $\pp \in X$. By the discussion above, $\pp$ is one of the three configurations with positions in $\{e, e+1\}$ (where $e = \eta(t)$), and its parent letter $w_t \in \{0,1,4\}$ satisfies $\varphi(w_t) = xy$. We compute $g(\pp) = (\cc, \uu, \vv)$ for each case:
		\begin{enumerate}
			\item If $\pp = (e, e, e, e+1)$, then $(d_0, d_1, d_2) = (\eps, \eps, x)$, which gives $\cc = (x, x, x, y)$, $\uu = \zero$, and $\vv = \psi(x)$.
			\item If $\pp = (e, e, e+1, e+1)$, then $(d_0, d_1, d_2) = (\eps, x, \eps)$, which gives $\cc = (x, x, y, y)$, $\uu = \psi(x)$, and $\vv = -\psi(x)$.
			\item If $\pp = (e, e+1, e+1, e+1)$, then $(d_0, d_1, d_2) = (x, \eps, \eps)$, which gives $\cc = (x, y, y, y)$, $\uu = -\psi(x)$, and $\vv = \zero$.
		\end{enumerate}
		Running over $xy \in \{03, 43, 01\}$ yields the nine stated triples.
	\end{proof}
	Explicitly, since $\psi(0) = \ee_0$ and $\psi(4) = \ee_4$, the nine starting states of $g(X)$ are:
	\begin{itemize}
		\item for $\varphi(0) = 03$:
		\[
		\bigl((0,0,0,3),\, \zero,\, \ee_0\bigr), \quad
		\bigl((0,0,3,3),\, \ee_0,\, -\ee_0\bigr), \quad
		\bigl((0,3,3,3),\, -\ee_0,\, \zero\bigr);
		\]
		\item for $\varphi(1) = 43$:
		\[
		\bigl((4,4,4,3),\, \zero,\, \ee_4\bigr), \quad
		\bigl((4,4,3,3),\, \ee_4,\, -\ee_4\bigr), \quad
		\bigl((4,3,3,3),\, -\ee_4,\, \zero\bigr);
		\]
		\item for $\varphi(4) = 01$:
		\[
		\bigl((0,0,0,1),\, \zero,\, \ee_0\bigr), \quad
		\bigl((0,0,1,1),\, \ee_0,\, -\ee_0\bigr), \quad
		\bigl((0,1,1,1),\, -\ee_0,\, \zero\bigr).
		\]
	\end{itemize}
	
	\begin{lemme}[Every additive cube ``descends'' from a starting state]\label{lem:descente}
		Let $\pp \in \N^4$ with $p_0 \le p_1 \le p_2 \le p_3$ and $p_0 < p_3$. There exists $k \ge 0$ such that $\parent^k(\pp) \in X$.
	\end{lemme}
	
	\begin{proof}
		From Lemma~\ref{lem:eta}(3) we know that the function $\parent$ is increasing (not necessarily strictly). Hence, we have $\parent^i(p_0) \le \parent^i(p_j) \le \parent^i(p_3)$ for every $j \in \{0,\dots,3\}$ and every $i \ge 0$.
		By Lemma~\ref{lem:eta}(3) again, $\parent^i(p_0) = \parent^i(p_3) = 0$ for sufficiently large $i$. Since $p_0 < p_3$, then there exists a maximal integer $k \ge 0$ such that $\parent^k(p_0) < \parent^k(p_3)$.
		By maximality of $k$, we have $\parent^{k+1}(p_0) = \parent^{k+1}(p_3)$.
		The inequalities then force all four coordinates of $\parent^k(\pp)$ to share the same parent, which proves that $\parent^k(\pp) \in X$.
	\end{proof}
	
	\subsection{Characterizing additive cubes by reachability in $\cG$}
	
	\begin{proposition}\label{prop:reduction}
		Under (H0):
		\begin{enumerate}
			\item[(a)] If $\xx$ contains an additive cube with boundaries $q_j = \pos(p_j,r_j)$, then $p_0 \le p_1 \le p_2 \le p_3$, $p_0 < p_3$, and there exists an integer $k \ge 0$ such that $\parent^k(\pp) \in X$. Moreover,
			\[
			g(\parent^k(\pp)),\, g(\parent^{k-1}(\pp)),\, \dots,\, g(\pp)
			\]
			is a path in $\cG$ from a vertex of $g(X)$ to a vertex of $Z_\tau$.
			\item[(b)] Conversely, if a path in $\cG$ leads from a vertex of $g(X)$ to a vertex of $Z_\tau$, then $\xx$ contains an additive cube.
		\end{enumerate}
	\end{proposition}
	
	\begin{proof}
		\begin{enumerate}
			\item[(a)] Since the positions $q_j$ are strictly increasing, one has $p_0 \le p_1 \le p_2 \le p_3$ (Lemma~\ref{lem:pos-tau}). If $p_0 = p_3$, then all four positions $p_j$ are equal, and since $r_3 < \ell(w_{p_0})$, the factor $\xx[q_0, q_3)$ is entirely contained in $\tau(w_{p_0})$, which contradicts (H0). Hence $p_0 < p_3$, and from Lemma~\ref{lem:descente} there exists an index $k \ge 0$ such that $\parent^k(\pp) \in X$. By Lemma~\ref{lem:automate-marche}(1), descending from $\parent^k(\pp)$ to $\pp$ in $\cT$ projects via $g$ onto a valid path in $\cG$. Finally, the equations~\eqref{eq:systeme-bords} characterizing the additive cube ensure that $g(\pp) \in Z_\tau$.
			
			\item[(b)] Conversely, consider now a path in $\cG$ from $g(\pp')$ (with $\pp' \in X$) to $(\cc,\uu,\vv) \in Z_\tau$. By Lemma~\ref{lem:automate-marche}(2), this path is the image under $g$ of a path of descendants in $\cT$ from $\pp'$ to some quadruple $\pp$, so that $g(\pp) = (\cc,\uu,\vv)$. Since $p'_0 < p'_3$ and $\eta$ is strictly increasing, the relative order is preserved and $p_0 < p_3$.
			
			Let $\rr$ be the offsets given by Definition~\ref{def:Ztau} and let $q_j = \pos(p_j, r_j)$. Since $c_j = w_{p_j}$, these are valid positions in $\xx$. The equations defining $Z_\tau$  imply $q_0 - 2q_1 + q_2 = 0$ and $q_1 - 2q_2 + q_3 = 0$, so the length $L := q_{j+1} - q_j$ is constant for $j \in \{0, 1, 2\}$. Furthermore,
			\[
			q_0 < \bl \cdot \sigma(p_0) + \ell(w_{p_0}) = \bl \cdot \sigma(p_0 + 1) \le \bl \cdot \sigma(p_3) \le q_3,
			\]
			which gives $3L = q_3 - q_0 > 0$, hence $L > 0$ and $q_0 < q_1 < q_2 < q_3$. The other two equations give the equality of the sums: $\xx[q_0,q_3)$ is an additive cube.
		\end{enumerate}
	\end{proof}
	
	Thus, searching for an additive cube in $\xx$ is equivalent to searching for a path in the \emph{infinite} graph $\cG$ from one of the nine vertices of $g(X)$ to a vertex of $Z_\tau$. It remains to make this graph finite.

	\section{Restricting the graph $\cG$ to a finite subgraph}\label{sec:fini}
	
	\subsection{Level zero and the transversality hypothesis (T)}
	
	\paragraph{Constraints at the leaf level.}
	If $\xx$ contains an additive cube, then the imbalance vectors $\uu,\vv$ and their sum $\uu + \vv = \psi(d_2) - \psi(d_0)$ (a Parikh difference between two factors of $\w$) satisfy
	\begin{equation}\label{eq:contrainte1}
		\Gamma \uu = -\bb(\cc,\rr), \qquad \Gamma \vv = -\bb'(\cc,\rr), \qquad \Gamma(\uu+\vv) = -\bigl(\bb(\cc,\rr) + \bb'(\cc,\rr)\bigr),
	\end{equation}
	and, by Theorem~\ref{thm:contraction},
	\begin{equation}\label{eq:contrainte2}
		|\theta_3(\uu)| \le C_3, \qquad |\theta_3(\vv)| \le C_3, \qquad |\theta_3(\uu + \vv)| \le C_3.
	\end{equation}
	Let $S \subset \Z^2$ be the finite set
	\[
	S = \left\{ -\bb(\cc,\rr), \; -\bb'(\cc,\rr), \; -\bigl(\bb(\cc,\rr) + \bb'(\cc,\rr)\bigr) \;\middle|\; \cc \in \Sigma^4, \; 0 \le r_j < \ell(c_j) \right\}.
	\]
	
	\begin{definition}[Initial state of level zero $F_0$]\label{def:S-F0}
		We denote 
		$F_0 = \left\{ \zz \in \Z^4 \;\middle|\; \Gamma \zz \in S \quad\text{and}\quad |\theta_3(\zz)| \le C_3 \right\}$. 
	\end{definition}
	If $\xx$ contains a cube, then $\uu, \vv, \uu + \vv \in F_0$.
	
	\paragraph{Why (T) is necessary: a cylinder cut by a hyperplane}
	For $F_0$ to be finite, an equation $\Gamma\zz = t$ ($t \in S$) and the bound $|\theta_3(\zz)| \le C_3$ must leave only finitely many integer points. For our eight morphisms, $\bs$ is not collinear with $\bl$: $\Gamma$ has rank $2$ and the set of real solutions of $\Gamma\zz = t$ is an affine plane $\zz_0 + \ker\Gamma$. This equation alone is not sufficient. On the other hand, the bound on $\theta_3$ places $\zz$ in the closed cylinder $\ker\theta_3 + K$, with axis $\ker\theta_3 = \operatorname{Vect}_\R(q_1,q_2) \subset \R^4$ and compact section $K = \operatorname{Vect}_\R(\operatorname{Re} q_3, \operatorname{Im} q_3) \cap \{|\theta_3| \le C_3\}$. The intersection of the affine plane and the cylinder is compact if and only if $\ker\Gamma \cap \ker\theta_3 = \{0\}$.
	
	\begin{definition}[Transversality hypothesis (T)]\label{def:T}
		The morphism 
		$\tau$ satisfies (T) if the real $4\times 4$ matrix
		\[
		\Phi = \begin{pmatrix} \Gamma \\ \operatorname{Re}(\theta_3) \\ \operatorname{Im}(\theta_3) \end{pmatrix}
		= \begin{pmatrix}
			\ell(0) & \ell(1) & \ell(3) & \ell(4) \\
			s(0) & s(1) & s(3) & s(4) \\
			\text{---} & \operatorname{Re}(\theta_3) & \text{---} \\
			\text{---} & \operatorname{Im}(\theta_3) & \text{---}
		\end{pmatrix}
		\]
		is invertible, that is, $\det(\Phi) \ne 0$ (equivalent to $\ker\Gamma \cap \ker\theta_3 = \{0\}$).
	\end{definition}
	
	\begin{exemple}[The matrix $\Phi$ for $\tau_1$]
		\[
		\Phi \approx \begin{pmatrix}
			2 & 3 & 1 & 2 \\[3pt]
			3 & 1 & 3 & 5 \\[3pt]
			0.5788 & -0.3219 & -0.0690 & -0.1662 \\[3pt]
			-0.5749 & 0.2183 & 0.6165 & -0.6810
		\end{pmatrix}, \qquad\det(\Phi) \approx -10.53 \ne 0.
		\]
	\end{exemple}
	
	\begin{lemme}[Finiteness of level zero]\label{lem:F0-fini}
		Under (T), the set $F_0$ is finite.
	\end{lemme}
	
	\begin{proof}
		$\Phi: \R^4 \to \R^4$ is a linear isomorphism. Identifying $\C$ with $\R^2$: $F_0 = \Z^4 \cap \Phi^{-1}\bigl(S \times \{\zeta \in \C: |\zeta| \le C_3\}\bigr)$, and the preimage of a compact set under an isomorphism is compact; its intersection with the lattice $\Z^4$ is therefore finite.
	\end{proof}
	
	\subsection{The box $\U$, the state space $H$ and the finite graph $\cG'$}
	Following the reduction strategy introduced in~\cite{CCSS}, we now bound all coordinates to restrict the exploration of $\cG$ to a finite induced subgraph. In~\cite{CCSS}, the initial bounds on $\theta_1$ and $\theta_2$ were deduced from a fixed two-dimensional integer lattice. Here, since $\tau$ and the boundary offsets modify lengths and sums, these initial bounds are instead controlled by the finite set $F_0$, whose finiteness is ensured by the transversality hypothesis~(T).
	
	\begin{definition}[Constants, box $\U$, space $H$ and graph $\cG'$]\label{def:boite-U}
		For $j \in \{1,2\}$, we set
		\[
		C_j = \max\left( \max_{\zz \in F_0} |\theta_j(\zz)|, \; \frac{\kappa_j}{|\lambda_j| - 1} \right), \qquad C_4 = C_3,
		\]
		then
		\[
		\U = \left\{ \zz \in \Z^4 \;\middle|\; |\theta_j(\zz)| \le C_j \text{ for } 1 \le j \le 4 \right\},
		\qquad
		H = \left\{ (\cc, \uu, \vv) \in \Sigma^4 \times \Z^4 \times \Z^4 \;\middle|\; \uu, \vv, \uu+\vv \in \U \right\}.
		\]
		$\cG'$ is the subgraph of $\cG$ induced on $H$, $A_0 = g(X) \cap H$, and $R$ is the set of vertices of $\cG'$ reachable from $A_0$.
	\end{definition}
	
	Each $\theta_j$ is bounded on $\U$ and $\zz = Q\theta(\zz)$ (Definition \ref{def:theta}), hence $\U$ is finite; $H$ and $\cG'$ are finite too, with at most $4^4 |\U|^2$ vertices.
	
	\begin{remarque}[One may always round up]\label{rem:mono}
		The rest of the work only uses the fact that $\U$ \emph{contains} the vectors that actually occur. One may therefore replace each $C_j$ by any upper bound: $\U$, $H$, $\cG'$, $R$ grow and Theorem~\ref{thm:reduction} remains valid. This is what the numerical computation does (Section~\ref{sec:calcul}).
	\end{remarque}
	
	\begin{proposition}[Confinement in $H$]\label{prop:confinement}
		Under the hypotheses of Proposition~\ref{prop:reduction}(a), all vertices $g(\parent^i(\pp))$, $0 \le i \le k$, belong to $H$.
	\end{proposition}
	
	\begin{proof}
		Set $d_j^{(i)} = \w[\parent^i(p_j), \parent^i(p_{j+1}))$ for $j = 0,1,2$. The three vectors $\uu^{(i)} = \psi(d_1^{(i)}) - \psi(d_0^{(i)})$, $\vv^{(i)} = \psi(d_2^{(i)}) - \psi(d_1^{(i)})$ and $\uu^{(i)} + \vv^{(i)} = \psi(d_2^{(i)}) - \psi(d_0^{(i)})$ are deviations $\Delta_i$ for a pair of factors.
		
		\emph{Level $0$.} By~\eqref{eq:contrainte1} and~\eqref{eq:contrainte2}, these three vectors are in $F_0$; hence $|\theta_j(\Delta_0)| \le \max_{F_0} |\theta_j|$ for $j = 1,2$.
		
		\emph{Levels $i \ge 0$.} For $j = 1,2$, Lemma~\ref{lem:remontee} (with $B = \max_{F_0}|\theta_j|$) gives $|\theta_j(\Delta_i)| \le C_j$; for $j = 3,4$, Theorem~\ref{thm:contraction} gives $|\theta_j(\Delta_i)| \le C_3$. The three vectors are therefore in $\U$.
	\end{proof}
	
	\begin{theoreme}[Reduction Theorem]\label{thm:reduction}
		Let $\tau: \Sigma^* \to \{0,1,2,3\}^*$ be non-erasing, satisfying (H0) and (T). Then
		\[
		\xx = \tau(\w) \text{ is without additive cubes} \iff R \cap Z_\tau = \emptyset.
		\]
	\end{theoreme}
	
	\begin{proof}
		($\Leftarrow$) Suppose that $\xx$ contains an additive cube. By Proposition~\ref{prop:reduction}(a), there exists a path in $\cG$ from a vertex of $g(X)$ to a vertex of $Z_\tau$. By Proposition~\ref{prop:confinement}, all its vertices are in $H$: it is a path of $\cG'$ issuing from $A_0$, and its endpoint is in $R \cap Z_\tau$.
		
		($\Rightarrow$) If $R \cap Z_\tau \ne \emptyset$, there exists a path of $\cG' \subset \cG$ from $A_0 \subset g(X)$ to $Z_\tau$; Proposition~\ref{prop:reduction}(b) provides an additive cube in $\xx$.
	\end{proof}
	
	Theorem~\ref{thm:reduction} shows that the exhaustive search is indeed a proof: if an additive cube exists in $\xx$, even arbitrarily far away, its genealogy will leave a trace in the finite graph $\cG'$, from one of the $9$ seeds to $Z_\tau$. For $\tau = \mathrm{id}$, we recover exactly the Theorem of \cite{CCSS}.

	\section{Computation and results}\label{sec:calcul}
	
	The \texttt{C++} program \texttt{proof.cpp} follows the proof:
	\begin{enumerate}
		\item it reads $\tau$ and tests (H0);
		\item it computes $M$, its eigenvalues, its eigenvectors of norm $1$, the forms $\theta_j$, and tests (T) ($\det\Phi \ne 0$);
		\item it computes $\mathcal{D}_9$ (in integers, via $\cQ$), $C_*$, $C_3$ and the $\kappa_j$;
		\item it enumerates $S$ and $F_0$, then computes $C_1$, $C_2$ and enumerates $\U$;
		\item it constructs $g(X)$ (Lemma~\ref{lem:gX});
		\item it performs a breadth-first traversal of $\cG'$ from $A_0$, the edges being generated ``on the fly'' by Definition~\ref{def:automate-G};
		\item for each vertex reached, it tests membership in $Z_\tau$ by running through the offsets $\rr$ (at most $81$ cases).
	\end{enumerate}
	Everything is exact (integers) except steps $2$ to $4$; the constants computed there are rounded upwards (relative margin $10^{-9}$), which is permitted by Remark~\ref{rem:mono}.
	
	For $\tau_1: 0 \mapsto 21,\ 1 \mapsto 100,\ 3 \mapsto 3,\ 4 \mapsto 23$:    
	\[
	\det(\Phi) \approx -10.53, \quad |S| = 43, \quad |F_0| = 67, \quad C_1 \approx 1.9032, 
	\]
    \[
	C_2 \approx 2.9818, \quad C_3 \approx 2.1758, \quad |\U| = 497.
	\]
    
The $9$ starting states are in $H$. The traversal reaches $|R| = 78340$ vertices, at maximal depth $23$, and $4078$ patterns $(\bb,\bb')$ are tested; no vertex reached is in $Z_\tau$. Theorem~\ref{thm:reduction} applies: $\tau_1(\w)$ contains no additive cube.

In fact, for all eight morphisms, the computation of $F_0$ gives $\max_{\zz \in F_0} |\theta_j(\zz)| \le \kappa_j / (|\lambda_j|-1)$ for $j \in \{1,2\}$. Consequently, $C_1$ and $C_2$ are entirely determined by the second term of the maximum in Definition~\ref{def:boite-U}, which depends only on $\varphi$. Hence the box $\U$, the state space $H$, the graph $\cG'$, and the set $R$ are identical for all eight morphisms; only $F_0$ and $Z_\tau$ change. The eight morphisms were obtained by an exhaustive search of candidates (images of length $\le 3$, prefix of $10^6$ letters without additive cubes):
	{\small
		\begin{align*}
			\tau_1 &: 0 \mapsto \phantom{0}21,\enspace 1 \mapsto 100,\enspace 3 \mapsto \phantom{00}3,\enspace 4 \mapsto \phantom{0}23, &
			\tau_2 &: 0 \mapsto \phantom{0}12,\enspace 1 \mapsto 233,\enspace 3 \mapsto \phantom{00}0,\enspace 4 \mapsto \phantom{0}10, \\[3pt]
			\tau_3 &: 0 \mapsto 122,\enspace 1 \mapsto 331,\enspace 3 \mapsto \phantom{00}0,\enspace 4 \mapsto \phantom{0}10, &
			\tau_4 &: 0 \mapsto 211,\enspace 1 \mapsto 002,\enspace 3 \mapsto \phantom{00}3,\enspace 4 \mapsto \phantom{0}23, \\[3pt]
			\tau_5 &: 0 \mapsto 021,\enspace 1 \mapsto 001,\enspace 3 \mapsto 323,\enspace 4 \mapsto 321, &
			\tau_6 &: 0 \mapsto 312,\enspace 1 \mapsto 332,\enspace 3 \mapsto 010,\enspace 4 \mapsto 012, \\[3pt]
			\tau_7 &: 0 \mapsto 132,\enspace 1 \mapsto 332,\enspace 3 \mapsto 101,\enspace 4 \mapsto 102, &
			\tau_8 &: 0 \mapsto 201,\enspace 1 \mapsto 001,\enspace 3 \mapsto 232,\enspace 4 \mapsto 231.
		\end{align*}
	}
	
	\begin{table}[htbp]
		\centering\small
		\begin{tabular}{lllcccc}
			\toprule
			$\tau(0), \tau(1), \tau(3), \tau(4)$ & $\bl$ & $\bs$ & $\det\Phi$ & $|F_0|$ & $|R|$ & $R \cap Z_\tau$ \\
			\midrule
			$21,\, 100,\, 3,\, 23$ & $2,3,1,2$ & $3,1,3,5$ & $-10.53$ & 67 & $78340$ & $\emptyset$ \\
			$12,\, 233,\, 0,\, 10$ & $2,3,1,2$ & $3,8,0,1$ & $+10.53$ & 67 & $78340$ & $\emptyset$ \\
			$122,\, 331,\, 0,\, 10$ & $3,3,1,2$ & $5,7,0,1$ & $+10.32$ & 119 & $78340$ & $\emptyset$ \\
			$211,\, 002,\, 3,\, 23$ & $3,3,1,2$ & $4,2,3,5$ & $-10.32$ & 119 & $78340$ & $\emptyset$ \\
			$021,\, 001,\, 323,\, 321$ & $3,3,3,3$ & $3,1,8,6$ & $-21.07$ & 89 & $78340$ & $\emptyset$ \\
			$312,\, 332,\, 010,\, 012$ & $3,3,3,3$ & $6,8,1,3$ & $+21.07$ & 89 & $78340$ & $\emptyset$ \\
			$132,\, 332,\, 101,\, 102$ & $3,3,3,3$ & $6,8,2,3$ & $+19.01$ & 87 & $78340$ & $\emptyset$ \\
			$201,\, 001,\, 232,\, 231$ & $3,3,3,3$ & $3,1,7,6$ & $-19.01$ & 87 & $78340$ & $\emptyset$ \\
			\bottomrule
		\end{tabular}
		\caption{Result of the traversal of $\cG'$ for the eight candidate morphisms.}
		\label{tab:resultats}
	\end{table}

    \section*{Disclaimer} 

ChatGPT (OpenAI, GPT-6) played a notable role in initiating this work. When consulted by the second author regarding the avoidability of additive cubes over $\{0,1,2,3\}$, it suggested a first candidate morphism, $\tau_1$, inspired by the morphism in~\cite{CCSS}, without a proof of correctness. The authors then proved that the infinite word $\tau_1(\w)$ is additive-cube-free and conducted an exhaustive search, which provided seven further morphisms as well as additional candidates of greater length, all of which can be proved by the same method. The authors take full responsibility for all mathematical results, proofs, and computations presented in this paper.

It is somewhat surprising that such an exhaustive search had not been done earlier; the emphasis in~\cite{Rao} on uniform morphisms of length at least~$16$ both discouraged any such search and directed intuition away from shorter and non-uniform ones, at least in the case of the first author.

	\bibliographystyle{plain}
	\bibliography{biblio}

@article{CCSS,
	author    = {Cassaigne, J. and Currie, J. D. and Schaeffer, L. and Shallit, J.},
	title     = {Avoiding three consecutive blocks of the same size and same sum},
	journal   = {Journal of the ACM},
	volume    = {61},
	number    = {2},
	pages     = {10:1--10:25},
	year      = {2014},
	note      = {arXiv:1106.5204}
}

@book{Lothaire,
	author    = {Lothaire, M.},
	title     = {Combinatorics on Words},
	series    = {Encyclopedia of Mathematics and its Applications},
	volume    = {17},
	publisher = {Cambridge University Press},
	address   = {Cambridge},
	year      = {1983}
}

@inproceedings{LR,
	author    = {Lietard, F. and Rosenfeld, M.},
	title     = {Avoidability of additive cubes over alphabets of four numbers},
	booktitle = {Developments in Language Theory (DLT 2020)},
	series    = {Lecture Notes in Computer Science},
	volume    = {12086},
	pages     = {192--206},
	publisher = {Springer},
	year      = {2020},
	doi       = {10.1007/978-3-030-48516-0_15}
}

@article{Rao,
	author    = {Rao, M.},
	title     = {On some generalizations of abelian power avoidability},
	journal   = {Theoretical Computer Science},
	volume    = {601},
	pages     = {39--46},
	year      = {2015}
}

@article{PV,
	author    = {Pirillo, G. and Varricchio, S.},
	title     = {On uniformly repetitive semigroups},
	journal   = {Semigroup Forum},
	volume    = {49},
	pages     = {125--129},
	year      = {1994}
}
	
\end{document}